
\documentstyle{amsppt}
\magnification=1200
\def\ss{\vskip.10in}
\def\ls{\vskip.25in}
\def\Hom{\Cal H\text{\rm om}}
\def\Ext{\text{\rm Ext}}

\topmatter
\title{\bf Unobstructedness of Calabi-Yau OrbiKleinfolds}
\endtitle
\def\C{\Bbb C}
\ls
\author{Z. Ran}
\endauthor
\affil  Department of Mathematics\\ University of California\\
Riverside, CA  92521, USA\\ ziv\@ucrmath.ucr.edu \endaffil
\thanks  ${}^*$ Supported in part by NSF DMS 9202050
\endthanks
\abstract{\it We show that Calabi-Yau spaces with certain types of
hypersurface-
quotient singularities have unobstructed deformations.  This applies in
particular to all Calabi-Yau orbifolds nonsingular in codimension 2.}
\endabstract
\endtopmatter
\document
\baselineskip=18pt

In recent years Calibi-Yau manifolds have attracted a great deal of
attention, motivated both by their role in Classification Theory and by
their connections with Physics, in particular the phenomenon of Mirror
Symmetry.   Both the Classification Theory and Physics viewpoints suggest
looking more generally at a suitable class of {\it singular} Calabi-Yau
spaces, e.g. most known Mirror Symmetry constructions involve forming a
quotient under a finite group action, thus leading naturally to Calabi-Yau
{\it orbifolds} and families of such.

Now a basic result about Calabi-Yau manifolds $X$ is the theorem of
Bogomolov-Tian-Todorov asserting that $X$ has unobstructed deformations.
Some extensions of this theorem to the case of singular $X$ have been
considered in [K], [R1], [T], [N].  However, these extensions do not
cover, e.g. the case of orbifolds.  Thus one is naturally led to ask, as
did D. Morrison, whether unobstructedness holds for Calabi-Yau orbifolds, say
nonsingular in codimension 2.   Our purpose here is to answer this question
affirmatively and, in fact, to prove a rather more general statement
which at the same time generalizes the result of [R1], allowing
(local) quotients of Kleinfolds, whence the term `orbiKleinfolds'.

The proof is particularly easy in the orbifold case, requiring little more
than a simple combination of (either form of) the `dual unobstructedness
criterion' of [R2] with the ideas surrounding Schlessinger's rigidity
theorem [S].  The more general orbiKleinfold case requires more work,
which we feel is justified by the interest of the result.

We begin by recalling that a {\it Kleimian} singularity is an
isolated simple hypersurface singularity.  These are classified by Arnold
into types A-D-E (cf. [1], p.~132).  We generalize this notion as follows.
\ss
\proclaim{Definition 1}  An isolated hypersurface singularity $(X, p)$ is
said to be {\rm weakly Kleinian} provided.
\roster
\item"(i)"  $X$ is $cDV$ (i.e., a general surface section through $p$ is a
DuVal or Kleinian singularity);

\item"(ii)"  $X$ admits a $\C^*$ action with positive weights;

\item"(iii)"  for a resolution $\tilde{X} \to X$, all components of the
exceptional locus are smooth divisors $E$ with $H^{p,q} (E) = 0$ for
$p\neq q$.
\endroster
\endproclaim

\demo{Remark}  Note that (i) implies that $X$ is a rational singularity, hence
$H^{p, 0} (E) = 0, \ p > 0$.  In particular, in dimension 3, (i) implies
(iii).
\enddemo

Next recall from [R1] that an isolated singularity $(X, p)$ is said to be
{\it good} if for a resolution $\tilde{X} \to X$ and $\Hat{\Tilde X}$ the
formal neighborhood of the exceptional locus, the map induced by exterior
derivative
$$
d: H^i (\Omega^j_{\Hat{\Tilde X}} ) \to H^i(\Omega^{j+1}_{\Hat{\Tilde X}})
$$
is injective in the range $i, j > 0, \ i \neq j, i + j < \dim \ X$.  All
3-fold rational singularities, and $A_1$ singularities in all dimensions are
good.  Presumably all Kleinian singularities are good, but this is
unproven.  The unobstructedness result of [R1] as stated allowed good
Kleinian singularities, it was remarked by Namikawa in [N] that the proof
covers the good weakly Kleinian case as well.
\ss
\proclaim{Definition 2}  An analytic variety $X$ is said to be a {\rm (weak)
OrbiKleinfold} if $X$ is locally of the form $V/G$, where $V$ has
(weak) Kleinian singularities and $G$ is a finite group acting
on $V$ which is `small' in that for all $g\in G, \ p \in V^g$, the induced
action of $g$ on the Zariski tangent $T_p V$ has no eigenspace of
codimension exactly 1.
\endproclaim
\demo{Remarks}
\roster
\item"1."  All orbifolds are orbiKleinfolds (cf. [St]).
\item"2."  In dimension 3, all terminal singularities are, by Mori's
classification (cf. [Rd]) (cyclic) quotients of isolated $cDV$ points,
hence `almost' weakly orbiKleinian; thus the class of weak orbiKleinfolds is
rather more general than that of orbifolds.
\item"3."  It seems reasonable that weak orbiKleinfolds are always
Cohen-Macaulay, but I cannot prove this.  In any case, the $CM$ property
certainly holds in all examples of interest, e.g. orbifolds and, more
generally, germs of the form $V(f)/G$ where $f$ is $G$-invariant:  indeed
$V(f)/G$ is then a Cartier divisor on an orbifold, hence $CM$.
\item"4."  If $X = V/G$ is a local weak orbiKleinfold, then $V$ admits
a $G$-equivariant retention $\tilde{V}$, yielding an orbifold
$X_{orb} = V/G$ with a birational morphism to $X$.  This construction
clearly globalizes, yielding an `orbifold resolution' $X_{orb} \to X$ for
any weak orbiKleinfold, which is an isomorphism off the (discrete)
non-orbifold locus of $X$.
\endroster
\enddemo
\ss
\proclaim{Definition 3}    {\rm A Calabi-Yau (weak) orbiKleinfold} is a compact
irreducible (weak) orbiKleinfold $X$ such that

(i)  $X$ admits a resolution of singularities by a K\"ahler manifold;

(ii)  $X$ is Cohen-Macauley;

(iii)  the dualizing sheaf $\omega_X \cong \Cal O_X $.

{\rm Our main result is the following}
\endproclaim
\ss
\proclaim{Theorem 4}  Any Calabi-Yau weak orbiKleinfold nonsingular in
codimension 2 has unobstructed deformations.
\endproclaim

\demo{Remark}  This generalizes Theorem 1 of [R1].  In dimension 3,
unobstructedness of Calabi-Yau spaces with terminal singularities has
been proven by Namikawa [N], but not all weak orbiKleinfolds are terminal,
nor conversely.
\enddemo

\demo{Proof of theorem}  Let $X$ be a Calabi-Yau weak orbiKleinfold of
dimension $n$, $\Theta = \Hom (\Omega_X^1, \Cal O_X)$ its tangent sheaf,
and $j: X_{\text{reg}} \to X$ the inclusion of the open subset of regular
points.  In view of our depth hypothesis, results of Schlessinger [S],
as exposed in ([A], I.9-10), show that we have an isomorphism of
deformation functors
$$
Def (X) \overset \sim\to\rightarrow \ Def\ (X_{\text{reg}}) .
$$
In particular the first-order deformation group
$$
T^1 (X) = \ Ext^1 (\Omega^1_X , \Cal O_X) \cong H^1 (X_{\text{reg}}, \Theta ) ,
$$
and, more importantly for us obstructions to deforming $X$ may be taken in the
group $H^2 (X_{\text{reg}}, \Theta)$.  The following is in essence due to
Schlessinger (and is valid more generally for quotients, nonsingular in
codimension 2, of isolated hypersurface singularities).
\enddemo

\proclaim{Lemma 5}  (i)  $R^1 {j_*} \Theta_{X_{\text{reg}}}$ is supported in
the nonorbifold locus of $X$;

(ii)  $R^2 {j_*} \Theta_{X_{\text{reg}}} = 0$.
\endproclaim

\demo{Proof}  We work locally on a small open set $U = V/G \subset X$, where
$V$ is a local hypersurface with one (possibly) singular point $p$.  Put
$$
\align
V_0 &= V\{p\}\\
V_{00} &= \ \text{subset of $V_0$ where $G$ acts freely.}
\endalign
$$
Note that $V/V_{00}$ has codimension $\ge 3$ in $X$ and $U_{\text{reg}} =
V_{00}/G$.  As in [A], we have
$$
H^i (V_{\text{reg}}, \Theta) = H^i (V_{00}, \Theta_{V_{00}})^G =
H^i (V_0, \Theta_{V_0})^G , \ \ i = 1,2,
$$
and $H^i(V_0, \Theta_{V_0})$ coincides with $T^i_{V,p}$ and in particular
vanishes if $p$ is regular or $i=2$. \qed
\enddemo

Now put
$$
\tilde{\Omega}_X = j_* \Omega_{X_{\text{reg}}} .
$$
This is a complex of reflexive sheaves on $X$ and by Lemma 2 we have
$$
H^2 (X_{\text{reg}}, \Theta_{X_{\text{reg}}}) = H^2 (X, j_*
\Theta_{X_{\text{reg}}}) \hookrightarrow \ \Ext^2 (\tilde{\Omega}^1_X, \Cal
O_X) ,
$$
hence the latter group is an obstruction group for deformations of $X$.
By Grothendieck duality, $\Ext^2 (\tilde{\Omega}_X^1, \Cal O_X) = \Ext^2
(\tilde{\Omega}^1_X, \omega_X)$ is dual to $H^{n-2} (\tilde{\Omega}_X^1)$.
Now the following result generalizes simultaneously results of
Steenbrink [St] and ([R1], Proposition 4).
\ss
\proclaim{Proposition 6}  Let $X$ be a compact weak orbiKleinfold
of dimension $n\ge 2, \ \pi_0: X_{\text{orb}} \to X$ an orbifold resolution
$\pi_1 : \tilde{X} \to X_{\text{orb}}$ a resolution of singularities with
$\tilde{X}$ K\"ahler.  Then

{\parindent=25pt
\item{(i)}  $\tilde{\Omega}^{\cdot}_{X} = \pi_{0*} \tilde{\Omega}^{\cdot}_{X_{
\text{orb}}} = (\pi_0 \circ \pi_1)_* \Omega^{\cdot}_{\tilde{X}} $;

\item{(ii)} $\tilde{\Omega}^{\cdot}_X$ is a resolution of the constant sheaf
$\C_X$;

\item{(iii)}  the Hodge-De Rham spectral sequence i
}
$$
E^{p,q}_1 = H^q (X, \tilde{\Omega}^p) \Rightarrow H^{p,q} (X, \C)
$$
degenerates at $E_1$ in degrees $\le n-1$.
\endproclaim

The argument that Proposition 6 implies Theorem 4, based on the
`$T^2$-injecting criterion' (Theorem 1.1, (ii) of [R2]) is identical to the
corresponding
argument in [R1].  In the case of orbifolds, Proposition 6 is due to Steenbrink
[St] (and consequently the reader only interested in orbifolds may stop
reading here).  Our proof of the Proposition combines ideas from [St] and
[R1].

\demo{Proof of Proposition 6}  (i)  The assertion is local, and at orbifold
points has been proven by Steenbrink ([St], Lemma 1.11).  For the
nonorbifold points, consider
$$
\gather
U = V/G \subset X , \\ \omega \in \Gamma (U, \tilde{\Omega}^i_X) =
\Gamma (U_{\text{reg}}, \Omega^i_{\text{reg}}) .
\endgather
$$
By ([St], Lemma 1.8), $\omega$ lifts to $\tilde{\omega} \in \Gamma
(V_{\text{reg}}, \Omega^i_{V_{\text{reg}}})$.  We may assume our
desingularization of $X$ is
locally over $U$ of the form $\tilde{V}/G$ where $\tilde{V} \to V$ is some
desingularization (which will also blow up the nonfree locus of $G$, in
addition to sing$V$).  By [R1], Proposition 4, (i), $\tilde{\omega}$
extends holomorphically to $\tilde{V}$ and being $G$-invariant descends to
$\tilde{V}/G$, as required.

(ii)  This is a consequence of ([St], Lemma 1.8), ([R1], Proposition 4), and
exactness of the functor of $G$-invariants.

(iii)  We prove the vanishing of the differentials $d_1^{ij}, \ i + j < n$,
the case of $d_r^{ij}, \ r\ge 2$, being similar.  Denote by
$\hat{\Omega}^j_X \subset \tilde{\Omega}^j_X$ the subsheaf of closed
forms and similarly for $\hat{\Omega}^j_{X_{\text{orb}}},
\hat{\Omega}^j_{\tilde{X}}$.  Consider the Leray spectral
sequence
$$
E^{p,q}_2 = H^p (X, R^q \pi_{0*} \tilde{\Omega}^j_{X_{\text{orb}}})
\Rightarrow H^i (X_{\text{orb}}, \tilde{\Omega}^j_{X_{\text{orb}}} )
$$
and the anlogous one for $\hat{\Omega}^j_{X_{\text{orb}}}$.  Using that
$\pi_0$ is an isomorphism off the finite nonorbifold locus of $X$, we get
a diagram
$$
\matrix
H^0 (R^{i-1} \pi_{0*} \tilde{\Omega}^j_{X_{\text{orb}}}) &\to&
H^0 (R^{i-1} \pi_{0*} \hat{\Omega}^{j+1}_{X_{\text{orb}}})&\overset
\beta_{ij}\to\rightarrow&H^0 (R^{i-1} \pi_{0*} \tilde{\Omega}^{j+1} ) \\
\downarrow&&\downarrow&&\downarrow\\
H^i (\tilde{\Omega}^j_X)&\to&H^i (\hat{\Omega}^{j+1}_X)&\to&
H^i (\tilde{\Omega}^{j+1}_X ) \\
\downarrow&&\downarrow&&\downarrow\\
H^i (\tilde{\Omega}^j_{X_{\text{orb}}}) &\overset \alpha_{ij}\to\rightarrow
&H^i(\hat{\Omega}^{j+1}_{X_{\text{orb}}})&\to&H^i(\tilde{\Omega}^{j+1}_{X_{
\text{orb}}})
\endmatrix
\tag*
$$
in which the composite of the middle-row arrow coincides with $d_{1}^{ij}$.
by Steenbrink [St], we have
$$
H^i (\hat{\Omega}^{j+1}_{X_{\text{orb}}}) = \Bbb H^i (F^{j+1}
\tilde{\Omega}^{\cdot}_{X_{\text{orb}}}) = F^{j+1} \Bbb H^{i+j+1}
(\tilde{\Omega}^{\cdot}_{X_{\text{orb}}}) ,
$$
where $F^{\cdot}$ denotes the Hodge (or stupid) filtration, and it follows
easily that $\alpha_{ij}$ vanishes.  To prove $d_1^{ij} = 0$ for
$i+j < n, \ j\neq i-2$, it will suffice, by $(*)$, to prove the  vanishing
of $\beta_{ij}$ in this range.  For this we use our goodness hypothesis.
Let $\hat{X}_{\text{orb}}$ be the formal neighborhood of the exceptional
locus of $\pi_0$ and similarly for $\Hat{\Tilde{X}}$.  Then we have a
diagram
$$
\matrix
H^{i-1} (\hat{X}_{\text{orb}}, \hat{\Omega}^{j+1}_{X_{\text{orb}}})&
\overset \hat{\beta}_{ij} \to \rightarrow&H^{i-1} (\hat{X}_{\text{orb}},
\Omega^{j+1}_{\text{orb}}) \\
\pi^*\quad\downarrow&&\pi^*\quad \downarrow\\
H^{i-1} (\Hat{\Tilde{X}}, \hat{\Omega}^{j+1}_{\tilde{X}})&\overset
\tilde{\beta}_{ij}\to \rightarrow&H^{i-1} (\Hat{\Tilde{X}},
\Omega^{j+1}_{\tilde{X}}) ,
\endmatrix
\tag**
$$
where $\hat{\beta}_{ij}, \tilde{\beta}_{ij}$ are the obvious maps.  By
goodness,
$\tilde{\beta}_{ij}$ vanishes whenever $i-1\neq j+1$, because the supposedly
injective exterior-derivative map vanishes on its image.  If we can prove the
vertical maps $\pi^*$ in $(**)$ are injective, it then follows that
$\hat{\beta}_{ij}$, hence $\beta_{ij}$, vanishes.  But the injectivity of
$\pi^*$ follows from Steenbrink's duality argument ([St], Proof of Thm. 1.12),
applied on the formal scheme $\hat{X}_{\text{orb}}$.
\ss
Finally, it remains to prove the vanishing of $d_1^{i,i-2}, 2 \le i \le
\frac{n+1}{2}$.  This is done just as in [R1], comparing $H^{\cdot} (X,\Bbb C)$
with $H^{\cdot} (X_{\text{orb}}, \Bbb C)$ and using the fact that the component
of the exceptional locus of $\pi_0$ are quotients of manifolds with
only $(p,p)$ cohomology, hence by Steenbrink's theory are themselves
orbifolds with all their cohomology of type $(p,p)$, hence contribute via
Gysin only $(p,p)$ classes to $H^{\cdot} (X_{\text{orb}})$. \qed
\enddemo
\ss
\demo{Remark}  After this was written, the author became aware of Fujiki's
paper
[F], which contains numerous results on symplectic orbifolds and their
deformations, but not including unobstructedness, which Fujiki
includes as an hypothesis in some statements.  In particular, in his
Theorem 4.8, p. 116, the hypothesis `$S$ is smooth' holds automatically.
\enddemo
\ss
\subheading{Acknowledgment}
\ss
We are grateful to Professor D.R. Morrison for suggesting the problem of
unobstructedness of Calabi-Yau orbifolds.
\ss
\subheading{References}
\roster
\ss
\item"[A]"  Artin, M.: `Deformations of singularities'.  Bombay:  Tata
Institute
1976.

\item"[F]"  Fujiki, A.:  `On primitively symplectic compact K\"ahler
$V$-manifolds
of dimension four' in:  {\it Classification of algebraic and analytic
manifolds}, K. Ueno, ed., Boston:  Birkhauser 1983, pp.~71-250.

\item"[K]"  Kawamata, Y.:  `Unobstructed deformations--a remark on a paper
by Z. Ran' (preprint).

\item"[L]"  Looijenga, E.J.N.:  `Isolated singular points on complete
intersections', Cambridge, Cambridge University Press, 1984.

\item"[N]"  Namikawa, Y.:  `On deformations of Calabi-Yau 3-folds with
terminal singularities', preprint.

\item"[R1]"  Ran, Z.:  `Deformations of Calabi-Yau Kleinfolds' in: S.T. Yau,
ed., `Essays on Mirrow Symmetry'.  Hong Kong, International Press 1992, pp.
451-457.

\item"[R2]"  $\underline{\text{\hskip.5in}}$:  `Hodge theory and
deformations of maps', Compositio Math. (to appear).

\item"[Rd]"  Reid, M.:  `A young person's guide to canonical singularities, in:
Proc. Symp. Pure Math, AMS.

\item"[S]"  Schlessinger, M.:  `Rigidity of quotient singularities'.
Invent.~Math. {\bf 14} (1971).

\item"[St]"  Steenbrink, J.H.M.:  `Mixed Hodge structure on the vanishing
cohomology' in:  `Real and complex singularities', Oslo 1976,
513--536.

\item"[T]"  Tian, G:  `Smoothing 3-folds with trivial canonical bundle
and ordinary double points',  ibid [R1] 458--479.
\endroster
\enddocument